\documentclass[12pt]{amsart}
\usepackage[latin1]{inputenc}
\usepackage{amscd}
\usepackage{latexsym}
\usepackage{pstcol,pst-plot,pst-3d}
\usepackage{multicol}
\usepackage{amsmath}
\usepackage{amssymb}
\usepackage{amsthm}
\psset{unit=0.7cm,linewidth=0.8pt,arrowsize=2.5pt 4}

\newpsstyle{fatline}{linewidth=1.5pt}
\newpsstyle{fyp}{fillstyle=solid,fillcolor=verylight}
\definecolor{verylight}{gray}{0.97}
\definecolor{light}{gray}{0.9}
\definecolor{medium}{gray}{0.85}

\def\NZQ{\Bbb}
\def\NN{{\NZQ N}}

\def\ZZ{{\NZQ Z}}

\def\CC{{\NZQ C}}

\def\FF{{\NZQ F}}
\def\GG{{\NZQ G}}

\def\frk{\frak}               

\def\mm{{\frk m}}

\def\Phi{{\frk n}}
\def\Phi{{\frk N}}

\def\opn#1#2{\def#1{\operatorname{#2}}} 
\opn\chara{char} \opn\length{\ell} \opn\pd{pd} \opn\rk{rk}
\opn\projdim{proj\,dim} \opn\injdim{inj\,dim} \opn\rank{rank}
\opn\depth{depth} \opn\grade{grade} \opn\height{height}
\opn\embdim{emb\,dim} \opn\codim{codim}

\opn\Tr{Tr} \opn\bigrank{big\,rank}
\opn\superheight{superheight}\opn\lcm{lcm}
\opn\trdeg{tr\,deg}
\opn\reg{reg} \opn\lreg{lreg} \opn\ini{in} \opn\lpd{lpd}
\opn\size{size} \opn\Pf{Pf} \opn\GL{GL} \opn\SL{SL} \opn\mod{mod}
\opn\ord{ord} \opn\Gin{Gin}
\opn\Hilb{Hilb}\opn\adeg{adeg}\opn\std{std}\opn\ip{infpt}
\opn\pol{pol}
\opn\div{div} \opn\Div{Div} \opn\cl{cl} \opn\Cl{Cl}
\opn\Spec{Spec} \opn\Supp{Supp} \opn\supp{supp} \opn\Sing{Sing}
\opn\Ass{Ass} \opn\Min{Min}
\opn\Ann{Ann} \opn\Rad{Rad} \opn\Soc{Soc}
\opn\Syz{Syz} \opn\Im{Im} \opn\Ker{Ker} \opn\Coker{Coker}
\opn\Am{Am} \opn\Hom{Hom} \opn\Tor{Tor} \opn\Ext{Ext}
\opn\End{End} \opn\Aut{Aut} \opn\id{id}

\opn\nat{nat}
\opn\pff{pf}
\opn\Pf{Pf} \opn\GL{GL} \opn\SL{SL} \opn\mod{mod} \opn\ord{ord}
\opn\Gin{Gin} \opn\Hilb{Hilb} \opn\cd{cd}
\opn\aff{aff} \opn\con{conv} \opn\relint{relint} \opn\st{st}
\opn\lk{lk} \opn\cn{cn} \opn\core{core} \opn\vol{vol}
\opn\link{link} \opn\star{star}\opn\limdepth{limdepth}
\opn\limdim{limdim}
\opn\gr{gr}

\def\pot#1#2{#1[\kern-0.28ex[#2]\kern-0.28ex]}

\opn\dirlim{\underrightarrow{\lim}}
\opn\inivlim{\underleftarrow{\lim}}
\let\tensor=\otimes
\let\iso=\cong

\let\Dirsum=\bigoplus

\let\to=\rightarrow

\def\Implies{\ifmmode\Longrightarrow \else
        \unskip${}\Longrightarrow{}$\ignorespaces\fi}
\def\implies{\ifmmode\Rightarrow \else
        \unskip${}\Rightarrow{}$\ignorespaces\fi}
\def\iff{\ifmmode\Longleftrightarrow \else
        \unskip${}\Longleftrightarrow{}$\ignorespaces\fi}

\let\:=\colon
\newtheorem{Theorem}{Theorem}[section]
\newtheorem{Lemma}[Theorem]{Lemma}
\newtheorem{Corollary}[Theorem]{Corollary}
\newtheorem{Proposition}[Theorem]{Proposition}
\newtheorem{Remark}[Theorem]{Remark}

\newtheorem{Definition}[Theorem]{Definition}

\def\Young#1{\vbox{\smallskip\offinterlineskip
    \halign{&\vbox{##}\kern-\Thickness\cr #1}}}

\newdimen\Squaresize \Squaresize=12pt
\newdimen\Thickness \Thickness=.3pt
\newdimen\Correction \Correction=7pt
\let\epsilon\varepsilon
\let\phi=\varphi
\let\kappa=\varkappa
\opn\dis{dis}
\def\pnt{{\raise0.5mm\hbox{\large\bf.}}}
\def\lpnt{{\hbox{\large\bf.}}}
\opn\Lex{Lex}

\begin{document}

\title{Local Duality for Bigraded Modules}

\author{J\"urgen Herzog and  Ahad Rahimi}
\address{J\"urgen Herzog, Fachbereich Mathematik und
Informatik, Universit\"at Duisburg-Essen, Campus Essen, 45117
Essen, Germany} \email{juergen.herzog@uni-essen.de}

\address{Ahad Rahimi, Fachbereich Mathematik und
Informatik, Universit\"at Duisburg-Essen, Campus Essen, 45117
Essen, Germany} \email{ahad.rahimi@uni-essen.de}

\address{}
 \maketitle

\begin{abstract}
In this paper we study local cohomology of finitely generated
bigraded modules over a standard bigraded ring with respect to the
 irrelevant  bigraded ideals and establish a duality theorem.
 Several applications are considered.
\end{abstract}
\bigskip

\section*{Introduction}
Let $R$ be a standard bigraded $K$-algebra with bigraded
irrelevant ideals $P$ generated by all elements of degree $(1,0)$,
and $Q$ generated by all elements of degree $(0,1)$. We want to
relate the local cohomology functors $H^i_P(-)$ and $H_Q^j(-)$ via
duality in the category of bigraded modules. In the ordinary local
duality theorem Matlis duality establishes isomorphisms between
the local cohomology modules of a module and its Ext-groups.

In our situation we have to consider Matlis duality for bigraded
modules. Given a  bigraded $R$-module $M$ we define the bigraded
Matlis-dual of $M$ to be  $M^\vee$  where the $(i,j)$th bigraded
component of $M^\vee$ is given by $\Hom_K(M_{(-i,-j)},K).$

As the main result of our paper we have the following duality
theorem:

\medskip
\noindent {\bf Theorem.} {\em Let $R$ be a standard bigraded
$K$-algebra
 with irrelevant
bigraded ideals $P$ and $Q$, and let $M$ be a finitely generated
bigraded $R$-module. Then there exists a convergent spectral
sequence
\[
E^2_{i,j}=H^{m-j}_P(H^i_{R_+}(M)^\vee )\underset{j} \Longrightarrow
H_Q^{i+j-m}(M)^\vee
\]
of bigraded $R$-modules, where $m$ is the minimal number of
homogeneous generators of $P$ and $R_+$ is the unique graded
maximal ideal of $R$. }

Note that the above spectral sequence degenerates when $M$ is
Cohen-Macaulay and one obtains  for all $k$ the following
isomorphims of bigraded $R$-modules
\begin{eqnarray}\label{1}
H^{k}_P(H^s_{R_+}(M)^\vee)\iso H^{s-k}_Q(M)^\vee \end{eqnarray}
where $s=\dim M$, see Corollary \ref{CM}.

Let $R_0$ be the $K$-subalgebra of $R$ which is generated by the
elements of bidegree $(1,0)$, and  let $N$ be any bigraded
$R$-module. Then for all $j$, the module $N_j=\Dirsum_iN_{(i,j)}$
is a  graded $R_0$-module with grading $(N_j)_i=N_{(i,j)}$.
Moreover if $N$ is finitely generated, then each $N_j$ is a
finitely generated $R_0$-module. In particular, if $M$ is an
$s$-dimensional Cohen-Macaulay module and if we set
$N=H^s_{R_+}(M)^\vee$, then $N$ is again an $s$-dimensional
Cohen-Macaulay module and by (\ref{1}) we obtain for all $j$ the
isomorphisms of graded $R_0$-modules
\begin{eqnarray}\label{2}
H^{k}_{P_0}(N_j)\iso (H^{s-k}_Q(M)_j)^\vee
\end{eqnarray}
 where $P_0$ is the
graded maximal ideal of $R_0$. Here we used, that $H^k_P(N)_j\iso
H^k_{P_0}(N_j)$ for all $k$ and $j$.

It has been conjectured by Brodmann and Hellus \cite{BH} that if $M$ is
finitely generated, then $H^{k}_Q(M)$ is tame, in other words, for
each $k$ there exists an integer $j_0$ such that either
$H^{k}_Q(M)_j=0$ for all $j\leq j_0$, or else $H^{k}_Q(M)_j\neq 0$ for
all $j\leq j_0$. In various cases this conjecture has been
confirmed, see \cite{BFL}, \cite{BH}, \cite{RS}, \cite{KS}, \cite{L} and \cite{Br} for a survey on this problem. In
case $M$ is Cohen-Macaulay the tameness problem translates, due to
(\ref{2}), to the following question: given a finitely generated
bigraded $R$-module $N$. Does there exist an integer $j_0$ such
that $H^{k}_{P_0}(N_j)=0$ for all $j\geq  j_0$, or else
$H^{k}_{P_0}(N_j)\neq 0$ for all $j\geq  j_0$?  More generally, one is
lead to conjecture that for a finitely generated graded
$R_0$-module $W$ and a finitely  generated bigraded $R$-module $N$
there exists for all $k$ an integer $j_0$ such that
$\Ext^k_{R_0}(N_j,W)=0$ for all $j\geq j_0$, or else
$\Ext^k_{R_0}(N_j,W)\neq 0$ for all $j\geq j_0$.

In Section 2 we use  our duality to give new proofs  of known
cases of the  tameness conjecture and also to add a few new cases
in which the conjecture holds, see Corollary \ref{tame1},
\ref{CMtame} and \ref{tame2}. The duality is also used in the
Corollaries \ref{structure1} and \ref{reg}   to prove some
algebraic properties of the modules $H^k_Q(M)_j$ in case $M$ is
Cohen-Macaulay.

\section{Proof the duality theorem}

Let $S=K[x_1, \dots,x_m, y_1,\dots,y_n]$ be the standard bigraded
polynomial ring over the field $K$. We set $K[x]=K[x_1,\dots,x_m]$
and $K[y]=K[y_1,\dots,y_n]$ and consider both as a standard graded
polynomial rings.

If $A$ is a standard (bi)graded $K$-algebra, and $M$ a (bi)graded
$A$-module.  Then we set $M^\vee=\Hom_K(M,K)$,
 and view  $M^\vee$
as (bi)graded $A$-module with the (bi)grading
\[
(M^\vee)_a=\Hom_K(M_{-a},K)
\]
for $a\in \ZZ$ (respectively $a\in \ZZ^2$ in the bigraded case).

The following simple fact is needed for the proof next lemma.

\begin{Lemma}
\label{tensor} Let $M$ be a graded $K[x]$-module and $N$ be a
graded $K[y]$-module. Then there exists a natural bigraded
isomorphism of bigraded $S$-modules
\[
(M \tensor_K N)^\vee \iso M^\vee \tensor_K N^\vee.
\]
\end{Lemma}
\begin{proof}
Let $S=K[x]\tensor_K K[y]=K[x,y]$. Note that  $M \tensor_K N$ is a
bigraded free $S$-module  with the natural bigrading
\[
(M \tensor_K N)_{(i,j)} =  M_i \tensor_K N_j.
\]
Thus we see that
\[
((M \tensor_K N)^\vee )_{(i,j)}=\Hom_K((M \tensor_K N)_{(-i,-j)},K)=\Hom_K(M_{-i} \tensor_K N_{-j},K).
\]
By using the universal property of tensor product one has the following natural isomorphism of $K$-vector spaces
\[
\Hom_K(M_{-i}\tensor_K N_{-j},K)\iso \Hom_K(M_{-i},K)\tensor_K \Hom_K(N_{-j},K).
\]
Thus we have
\begin{eqnarray*}
((M \tensor_K N)^\vee )_{(i,j)}&\iso &\Hom_K(M_{-i},K)\tensor_K \Hom_K(N_{-j},K)\\
                               &=&   (M^\vee)_i \tensor_K (N^\vee)_j\\
                               &=&   (M^\vee \tensor_K N^\vee)_{(i,j)}.
\end{eqnarray*}
So the desired isomorphism follows.
\end{proof}

\begin{Lemma}
\label{simple} Let $S=K[x_1, \dots,x_m, y_1,\dots,y_n]$ be the
standard bigraded polynomial ring over the field $K$ with the
irrelevant bigraded ideals $P=(x_1, \dots,x_m)$ and
$Q=(y_1,\ldots, y_n)$. Then we have the following isomorphism of
bigraded  $S$-modules
\[
H^m_P(\omega_S) \iso H^n_Q(S)^\vee ,
\]
where $\omega_S$ is the bigraded canonical module of $S$.
\end{Lemma}
\begin{proof}
First we notice that there is a  natural isomorphism of bigraded
$S$-modules
\[
H^m_P(S)\iso H^m_P(K[x])\tensor _K K[y].
\]
By the graded version of the local duality theorem  (see \cite
[Example 13.4.6]{BS}) we have
\[
H^m_P(K[x])^\vee \iso K[x](-m).
\]
Thus we see that
\begin{eqnarray*}
H^m_P(\omega_S)=H^m_P(S(-m,-n))&=&H^m_P(S)(-m,-n)\\
                               &\iso& (K[x](-m)^\vee \tensor_K K[y])(-m,-n)\\
                               &=& K[x]^\vee \tensor_K K[y](-n).
\end{eqnarray*}
On the other hand, using again the local duality theorem, Lemma
\ref {tensor} yields
\begin{eqnarray*}
H^n_Q(S)^\vee \iso (K[x]\tensor_K H^n_Q(K[y]))^\vee &\iso& (K[x]\tensor_K K[y](-n)^\vee)^\vee\\
                                                    &\iso& K[x]^\vee \tensor_K K[y](-n)^{\vee\vee}\\
                                                    &\iso& K[x]^\vee \tensor_K K[y](-n),
\end{eqnarray*}

as desired.
\end{proof}

\begin{Corollary}
\label{free}

Let $F$ by a finitely generated  bigraded free $S$-module, and set
$F^*=\Hom_S(F, \omega_S)$. Then there exists a natural isomorphism
of bigraded $S$-modules,
\[
H^m_P(F^*) \iso H^n_Q(F)^\vee.
\]
\end{Corollary}

\begin{proof}
Let $F=\Dirsum_{k=1}^{t} S(-a_{k},
-b_{k})$. Thus  $F^*=\Dirsum_{k=1}^{t}(\omega_S)(a_{k},b_{k})$ and hence by Lemma \ref{simple} we have
\begin{eqnarray*}
H^m_P(F^*)\iso  \Dirsum_{k=1}^{t} H^m_P(\omega_S)(a_{k},b_{k})& \iso & \Dirsum_{k=1}^{t} H^n_Q(S)^\vee(a_{k},b_{k})\\
&\iso & H_Q^n(\Dirsum_{k=1}^{t} S(-a_{k}, -b_{k}))^\vee \\ &\iso& H_Q^n(F)^\vee.
\end{eqnarray*}
\end{proof}

The previous result  can easily be extended as follows

\begin{Lemma}
\label{irrelevant}  Let  $\FF$ be a bounded complex of bigraded
free $S$-modules. We set $\FF^*=\Hom_S(\FF, \omega_S)$. Then we
have a functorial isomorphism
\[
H^m_P(\FF^*) \iso H^n_Q(\FF)^\vee
\]
of complexes of bigraded modules.
\end{Lemma}
\begin{proof}
In order to prove that the complexes of
$H_P^m(\FF^*)$ and $H_Q^n(\FF)^\vee$ are isomorphic,  we observe  that for any
bihomogeneous linear map $\phi:G \rightarrow F$ between finitely
generated free bigraded $S$-modules we obtain the following
commutative diagram
\[
\begin{array}{ccc}
 H_Q^n(F)^\vee & \stackrel  {\psi^\vee_1} \longrightarrow  &  H_Q^n(G)^\vee \\
\downarrow  &  & \downarrow  \\
H_P^m(F^*)        &   \stackrel {\psi_2}\longrightarrow &
H_P^m(G^*),
\end{array}
\]
where   $\psi_1=H_P^n(\phi)$ and $\psi_2=H_Q^m(\phi^*) $ and where the
vertical maps are  the isomorphisms given in Corollary \ref{free}.
The commutativity of the diagram results from the fact that all
maps in the diagram are functorial.
\end{proof}

\begin{Proposition}
\label{Duality} Let $M$ be a finitely generated bigraded
$S$-module, $P$ and $Q$ be the irrelevant bigraded ideals of $S$.
Then we have the following convergent spectral sequence
\[
E_{i,j}^2=
H^{m-j}_P(\Ext_S^{n+m-i}(M,\omega_S))\underset{j}\Longrightarrow
H_Q^{i+j-m}(M)^\vee.
\]
\end{Proposition}

\begin{proof}
Let $(\FF,d)$ be a bigraded free resolution of $M$ of length
$n+m$, and let $\GG$ be the complex of bigraded $S$-modules with
$G_i=\Hom_S(F_{m+n-i},\omega_S)$ and  differential
$\partial_i=\Hom_S(d_{m+n-i},\omega_S)$. Next we choose a bigraded
free resolution $\CC$ of the complex $\GG$. In other words, $\CC$
is a double complex $C_{ij}$ of finitely generated bigraded free
$S$-modules with $i,j\geq 0$ such that:
\begin{enumerate}
\item[(i)] the $i$th column of $\CC$ is a free resolution of $G_i$
for all $i$, i.e.

\[
\ H_j(C_{i\lpnt})=\left\{
\begin{array}{ll}
G_i  & \text{for $j=0$, }\\
 0  & \text{for  $j>0$.}
\end{array}
\right.
\]
\item[(ii)] for each row the image of $C_{i-1, j}\longleftarrow
C_{i,j}$ is a bigraded free direct summand of the kernel of
$C_{i-2,j}\longleftarrow C_{i-1,j}$. In particular, the homology
of
\[
C_{i-2,j}\longleftarrow C_{i-1,j}\longleftarrow
C_{i,j}
\]
 is a bigraded free $S$-module for all $i$ and $j$.

\item[(iii)] for each $i$ the complex
\[
0\longleftarrow H_i(C_{\lpnt,0}) \longleftarrow H_i(C_{\lpnt,1})
\longleftarrow  H_i(C_{\lpnt,2}) \longleftarrow\cdots
\]
is a bigraded free resolution of $H_i(\GG)$.
\end{enumerate}

Now we compute the total homology of the double complex $H_P^m(\CC)$:
Since all $G_i$ are free $S$-modules, it follows that the
complexes
\[
0 \longleftarrow G_i \longleftarrow C_{i,0}\longleftarrow
C_{i,1}\longleftarrow \cdots
\]
are all split exact. Hence the complexes
\[
0 \longleftarrow H^m_P(G_i) \longleftarrow H^m_P(C_{i,0})
\longleftarrow   H^{m}_P(C_{i,1}) \longleftarrow  \cdots
\]
are  again exact.

This implies that the $E^1$-terms of the double complex
$H_P^m(\CC)$ with respect to the column filtration are

\[
\ E^1_{i,j} =\left\{
\begin{array}{ll}
H^m_P(G_i)  & \text{for $j=0$,}\\
 0 & \text{for  $j>0$.}
\end{array}
\right.
\]
As a consequence, for the  $E^2$-terms  of $H_P^m(\CC)$ we have
that $E^2_{i,j}=0$ for $j>0$, and that $E^2_{i,0}$ is the $i$th
homology of the complex $H^m_P(\GG)$. Now we use Lemma
\ref{irrelevant} as well as \cite[Theorem, 1.1 ]{AR} and obtain
\[
\ E^2_{i,j} =\left\{
\begin{array}{ll}
H^{i-m}_Q(M)^\vee & \text{for $j=0$,}\\
 0  & \text{for  $j>0$,}
\end{array}
\right.
\]
since $H_i(H_P^m(\GG))=(H_{n+m-i}(H^n_Q(\FF)))^\vee$. From this it
follows that the $(i+j)$th total homology of $H_P^m(\CC)$ is equal
to $H^{i+j-m}_Q(M)^\vee$.

Now we compute the homology of $H_P^m(\CC)$ using the row
filtration. Each row $H^m_P(C_{\lpnt j})$ of $H_P^m(\CC)$ is split
exact with homology $H_i(H^m_P(C_{\lpnt j}))=H^m_P(H_i(C_{\lpnt
j}))$. In other words, $E^1_{i,j}=H^m_P(H_i(C_{\lpnt j}))$. Hence
by  property (iii) of the complex $\CC$  and by \cite[Theorem, 1.1
]{AR} it follows that
$E^2_{i,j}=H_P^{m-j}(\Ext^{m+n-i}_S(M,\omega_S))$. This yields the
desired conclusion.
\end{proof}

Now our main theorem is an easy consequence  of Proposition
\ref{Duality}:
\begin{proof}
As $R$ is a standard bigraded $K$-algebra,  it is the homomorphic
image of a standard bigraded polynomial ring
$S=K[x_1,\dots,x_m,y_1,\dots,y_n]$.  We may consider $R$ and $S$
as well as a standard graded $K$-algebras with the unique graded
maximal ideal $R_+$ (resp.\ $S_+$), and $M$ as a graded $R$-module
(resp.\ $S$-module). Then by the graded local duality theorem we
have
\[
 \Ext_S^{m+n-i}(M,\omega_S)\iso H^i_{S_+}(M)^\vee.
\]
Since $H^i_{S_+}(M)\iso H^i_{R_+}(M)$, it follows that
\[
H^{m-j}_P(H^i_{R_+}(M)^\vee )=H^{m-j}_{P}(\Ext_S^{m+n-i}(M,\omega_S)).
\]
Let $(x)=(x_1,\dots,x_m)$ and $(y)=(y_1,\dots,y_n)$ be the
irrelevant ideals of $S$. We note that
$H^{m-j}_{P}(\Ext_S^{m+n-i}(M,\omega_S))=H^{m-j}_{(x)}(\Ext_S^{m+n-i}(M,\omega_S))$
and that $ H_Q^{i+j-m}(M)^\vee=H_{(y)}^{i+j-m}(M)^\vee. $
Therefore, Proposition \ref{Duality} yields the desired convergent
spectral sequence.
\end{proof}

\begin{Corollary} {\em Let $R$ be a standard bigraded
$d$-dimensional Cohen-Macaulay $K$-algebra  with irrelevant
bigraded ideals $P$ and $Q$, and let $M$ be a finitely generated
bigraded $R$-module. Then there exists a convergent spectral
sequence
\[
E^2_{i,j}=H^{m-j}_P(\Ext_R^{d-i}(M,\omega_R))\underset{j} \Longrightarrow
H_Q^{i+j-m}(M)^\vee
\]
of bigraded $R$-modules, where $m$ is the minimal number of
homogeneous generators of $P$. }
\end{Corollary}
\begin{proof}
The assertion follows from our main theorem by using the fact that $H^i_{R_+}(M)^\vee=\Ext_R^{d-i}(M,\omega_R)$.
\end{proof}
\section{Some applications}

In this section, unless otherwise stated,  $R$ denotes a standard
bigraded  $K$-algebra of dimension $d$,  and $M$ a finitely
generated and bigraded $R$-module.

We note that for the $E^2$-terms in the  spectral sequence of our
main theorem we have $E^2_{i,j}=H^{m-j}_P(H^i_{R_+}(M)^\vee)=0$ if
$i<\depth M$ or $i>\dim M$ or $j<0$ or $j>m$. Thus the possible
non-zero $E^2$-terms are in the shadowed region of the following
picture.

 \begin{center}
\psset{unit=1.0cm}
\begin{pspicture}(0,0)(6,6)
 \psline(0.5,1)(0.5,5)
 \psline(0,1.5)(5.5,1.5)
 \pspolygon[style=fyp, fillcolor=light](1.5,1.5)(1.5,3.5)(4.5,3.5)(4.5,1.5)
 \psline[linestyle=dashed](1.5,1.5)(1.5,3.5)
 \psline[linestyle=dashed](1.5,3.5)(4.5,3.5)
 \psline[linestyle=dashed](4.5,3.5)(4.5,1.5)
 \psline[linestyle=dashed](0.5,2.5)(3,2.5)
 \psline[linestyle=dashed](3,2.5)(3,1.5)
 \psline[linestyle=dashed](0.0,3.5)(0.5,3)
 \psline[linestyle=dashed](0.5,3)(1,2.5)
 \psline[linestyle=dashed](1,2.5)(1.5,2)
 \psline[linestyle=dashed](1.5,2)(2,1.5)
 \psline[linestyle=dashed](2,1.5)(2.5,1)
 \rput(0.3,2.5){\blue{$j$}}
 \rput(3,1.2){\blue{$i$}}
 \rput(1.5, 1.5){\red{$\bullet$}}
 \rput(4.5, 3.5){\red{$\bullet$}}
 \rput(4.5, 3.9){{$(s,m)$}}
 \rput(3, 2.5){\red{$\bullet$}}
 \rput(1.5, 1.1){$(t,0)$}
 \rput(4.5,1.2){$s$}
 \rput(3.5, 2.5){\blue{$E^2_{i,j}$}}
 \rput(3,0.5){Figure 1}
 \end{pspicture}
\end{center}

Here $t=\depth M$ ,  $s=\dim M$ and
$E^2_{i,j}=H^{m-j}_P(H_{R_+}^i(M)^\vee)$.

We first observe that the graded local duality theorem is a
special case of our main theorem. In fact, if we assume that
$P=(0)$, then $m=0$, and  $\mm= Q$ is the unique graded maximal
ideal of $R$. Moreover,
 $E^2_{i,j}=E^\infty_{i,j}=0$ for $j\neq 0$ and  all $i$, since
 $H_{(0)}^k(-)=0$ if $k\neq 0$.
Therefore we have
\[
\Ext_R^{d-i}(M,\omega_R)=H^0_{(0)}(\Ext_R^{d-i}(M,\omega_R))\iso
H^i_\mm(M)^\vee.
\]

Considering Figure 1 we immediately obtain the following corner
isomorphisms

\begin{Proposition}
\label{corner} Let  $\dim M=s$ and $\depth M=t$. Then there are
natural isomorphisms
\[
H^{m}_P(H_{R_+}^t(M)^\vee)\iso H^{t-m}_Q(M)^\vee \quad\text{and}
\quad H^{0}_P(H_{R_+}^s(M)^\vee)\iso H^{s}_Q(M)^\vee.
\]
Moreover  for $i<t-m$ we have $H^{i}_Q(M)=0$.
\end{Proposition}

\begin{Definition}
\label{deftame}
{\em Let $R_0$ be a commutative Noetherian ring, $R$ a graded $R_0$-algebra and $N$ a graded $R$-module. The $R$-module $N$ is called {\em tame}, if there exists an integer $j_0$ such that either
\[
N_j=0\quad\text{for all}\quad j\leq j_0, \quad \text{or}\quad  N_j\neq 0\quad \text{for all}\quad j\leq j_0.
\]}
\end{Definition}

In  case of a standard bigraded $K$-algebra $R$ we let $R_0$ be
the $K$-subalgebra of $R$ generated by all elements of degree
$(1,0)$. Then $R$ is a graded $R_0$-algebra with  components
$R_j=R_{(*,j)}=\Dirsum_iR_{(i,j)}$. Let $N$ be a bigraded
$R$-module. We may view $N$ as a graded $R$-module with graded
components $N_j= N_{(*,j)}=\Dirsum_iN_{(i,j)}$. Each of the
modules $N_j$ is a graded $R_0$-module, and if $N$ is a finitely
generated $R$-module then each $N_j$ is a finitely generated
$R_0$-module.

Now let $M$ be a finitely generated bigraded $R$-module. Then
$H^i_P(M)_j=H^i_{P_0}(M_j)$ where $P_0$ is the graded maximal
ideal of $R_0$. Since $M_j$ is a finitely generated $R_0$-module
it follows that $H^i_{P_0}(M_j)$ is a graded Artinian
$R_0$-module. Hence we see that
$(H^i_P(M)^\vee)_j=(H^i_P(M)_{-j})^\vee=H^i_{P_0}(M_{-j})^\vee$ is
a finitely generated graded $R_0$-module for all $j$. Of course
this does not imply that  $H^i_P(M)^\vee$ is a finitely generated
$R$-module.

We denote by $\cd(M)$ the cohomological dimension of $M$ with
respect to $Q$, i.e.\ the number
\[
\cd(M)=\sup \{i\in \NN_0: H^i_Q(M)\neq 0 \}.
\]
\begin{Corollary}
Let $N=H_{R_+}^s(M)^\vee$. Then the following statements hold:
\begin{enumerate}
\item [(a)] $\cd(M)<\dim(M)$ if and only if $\depth_{R_0}N_j>0$
for all $j$; \item [(b)] if  $\cd(M)<\dim (M)-1$, then
$\depth_{R_0}N_j>1$ for all $j$.
\end{enumerate}
\end{Corollary}
\begin{proof}
We note that  $\cd(M)<\dim(M)$, if and only if  $H^{s}_Q(M)=0$.
Hence Proposition \ref{corner} yields part  (a) of the corollary.

For the proof of (b) we notice  that $H^{i}_P(N)=E_{s,m-i}^\infty$
for $i=0,1$, and that $E_{s,m-i}^\infty$ is a submodule of
$H_Q^{s-i}(M)^\vee$ for all $i$.  Thus our assumption implies that
$H^{i}_{P_0}(N_j)=H^{i}_{P}(N)_j=0$ for $i=0,1$ and all $j$. This
yields the desired conclusion.
\end{proof}

The second statement of the next corollary of  is well-known (see
\cite[Theorem, 4.8 (e)]{Br}).
\begin{Corollary}
\label{tame1} Let $M$ be a finitely generated bigraded $R$-module
of dimension $s$ and depth $t$. Then   $H^{t-m}_Q(M)$ and
$H^{s}_Q(M)$ are tame.
\end{Corollary}

\begin{proof}
We first prove $H^{t-m}_Q(M)$ is tame.  We set
$N=H_{R_+}^t(M)^\vee$ and $s_0=\dim_{R_0}N_j$ for $j\gg 0.$ Note
that $s_0\leq \dim R_0\leq m.$ Thus we have
$H^m_P(N)_j=H^m_{P_0}(N_j)=0$ for $j\gg 0$ if $s_0<m $ and
$H^m_P(N)_j=H^m_{P_0}(N_j)\neq 0$ for $j\gg 0$ if $s_0=m.$
Therefore by Proposition \ref {corner} there exists an integer
$j_0$ such that either
\[
H^{t-m}_Q(M)_j=0\quad\text{for all}\quad j\leq j_0,\quad \text{or}\quad H^{t-m}_Q(M)_j\neq 0\quad \text{for all}\quad j\leq j_0,
\]
as desired. In order to prove that $H^{s}_Q(M)$ is tame,  we set
$N=H_{R_+}^s(M)^\vee$.  Since  $H_{R_+}^s(M)$ is a graded Artinian
$R$-module, $N$ is a finitely generated graded $R$-module. Thus
$N_j$ is a finitely generated $R_0$-module. By \cite[Proposition,
2.5]{B} the set of associated prime ideals of $\Ass_{R_0}(N_j)$
is constant for large $j$. If  $P_0\in \Ass_{R_0}(N_j)$ it follows
that $H^{0}_P(N)_j=H^{0}_{P_0}(N_j)\neq 0$ for large $j$, and if
$P_0\notin \Ass_{R_0}(N_j)$ then $H^{0}_P(N)_j=H^{0}_{P_0}(N_j)=0$
for large $j$. Thus in view of Proposition \ref {corner},
$H^{s}_Q(M)$ is also tame.
\end{proof}

We say that $M$ is a generalized Cohen-Macaulay $R$-module if
$H^i_{R_+}(M)$ has finite length for all $i\neq \dim M$.

\begin{Proposition}
\label{generalized}  Let  $M$ be a generalized Cohen-Macaulay
$R$-module of dimension $s$. Then we have the following long exact
sequence of bigraded $R$-modules

\begin{eqnarray*}
0 \rightarrow H^{1}_P(H_{R_+}^s(M)^\vee) \rightarrow H^{s-1}_Q(M)^\vee \rightarrow H_{R_+}^{s-1}(M)^\vee \rightarrow\\
H^{2}_P(H_{R_+}^s(M)^\vee)\rightarrow
H^{s-2}_Q(M)^\vee \rightarrow H_{R_+}^{s-2}(M)^\vee \rightarrow\\
 \dots \rightarrow H^{s-m}_Q(M)^\vee\rightarrow H_{R_+}^{s-m}(M)^\vee\rightarrow
 0.
\end{eqnarray*}
 Moreover, we have the following isomorphisms
\[
H^i_{R_+}(M)\iso H^{i}_Q(M) \quad\text{for all $i<s-m$.}
\]
\end{Proposition}
\begin{proof}
Since $M$ is a generalized Cohen-Macaulay module, we have that
$H_{R_+}^i(M)^\vee$ is of finite length for $i\neq s.$ Thus
by Grothendieck's vanishing theorem \cite[Theorem, 6.1.2]{BS} we
see that $E^2_{i,j}=E^\infty_{i,j}=0$ for $j=0,\dots, m-1$ and
$i\neq s.$ The following picture will make this clear.

 \begin{center}
\psset{unit=1.0cm}
\begin{pspicture}(0,0)(6,6)
 \psline(0.5,1)(0.5,5)
 \psline(0,1.5)(6,1.5)

 \psline(1.5,3.5)(4.5,3.5)
 \psline(4.5,3.5)(4.5,1.5)

 \psline[linestyle=dashed](3.5,3.5)(4.5,2.5)
 \psline[linestyle=dashed](4.5,2.5)(5.5,1.5)
 \psline[linestyle=dashed](3.5,3.5)(3.5,1.5)
 \psline[linestyle=dashed](3,4)(3.5,3.5)
 \psline[linestyle=dashed](5.5,1.5)(6,1)
 \rput(3.5, 3.5){\red{$\bullet$}}
\rput(5, 2.5){{$(s,r)$}}
\rput(4, 3.8){{$(l,m)$}}
 \rput(4.5, 2.5){\red{$\bullet$}}
 \rput(5.5, 1.2){{$k$}}
 \rput(1.5, 1.2){$t$}
 \rput(4.5,1.2){$s$}
 \rput(3.5,1.2){$l$}
 \rput(3,0.5){Figure 2}
 \end{pspicture}
\end{center}
Therefore for all $k$ with $s\leq k<s+m$
 we get the following exact sequences
\[
0 \rightarrow E^\infty_{s,r}\rightarrow H^{l}_Q(M)^\vee \rightarrow E^\infty_{l,m}\rightarrow 0,
\]

\[
0 \rightarrow E^\infty_{l,m}\rightarrow E^2_{l,m} \stackrel
{}\rightarrow  E^2_{s,r-1}\rightarrow
E^\infty_{s,r-1}\rightarrow 0,
\]
where $l$ and $r$ are defined by the equations $s+r=l+m=k$.

Composing these two exact sequences we get the long exact sequence
\begin{eqnarray*}
\dots \rightarrow E^2_{s,r} \rightarrow H^{l}_Q(M)^\vee \rightarrow E^2_{l,m}\rightarrow E^2_{s,r-1}\rightarrow \\ H^{l-1}_Q(M)^\vee \rightarrow E^2_{l-1,m}\rightarrow  E^2_{s,r-2}\rightarrow \dots,
\end{eqnarray*}
which yields the desired exact sequence, observing that
$$H^0_P(H_{R_+}^i(M)^\vee)=H_{R_+}^i(M)^\vee$$
for $i\neq s$   since for such $i$ the modules
$H_{R_+}^i(M)^\vee$ have finite length. The last statement
of the proposition follows similarly.
\end{proof}

\begin{Corollary}
 \label{CM} Suppose $M$ is a generalized Cohen-Macaulay module of
dimension $s$. Then the  following conditions are equivalent:
\begin{enumerate}
\item[(a)] $M$ is Cohen-Macaulay;

\item[(b)] $H^{k}_P(H^s_{R_+}(M)^\vee)\iso H^{s-k}_Q(M)^\vee$ for
all $k$.
\end{enumerate}
\end{Corollary}

\begin{proof} (a)\implies (b): Since $M$ is Cohen-Macaulay we have
$H_{R_+}^i(M)^\vee=0$ for all $i\neq s$. Therefore it
follows from the long exact sequence in Proposition
\ref{generalized} that $H^{k}_P(H_{R_+}^s(M)^\vee)\iso
H^{s-k}_Q(M)^\vee$ for $k=1,\dots,m$. The assertion for $k=0$
follows from Proposition \ref{corner}.  The assertion is also
clear when $k<0$. Now assume that $k>m$. Then $s-k<s-m$, and hence
by Proposition \ref{generalized} it follows that
$H^{s-k}_Q(M)^\vee=H_{R_+}^{s-k}(M)^\vee=0$. On the other hand, we
also have $H^{k}_P(H_{R_+}^s(M)^\vee)=0$ because $k>m$.

(b)\implies (a): is proved the same way.
\end{proof}

As a generalization of Lemma \ref {simple} we obtain as an
immediate consequence of Corollary  \ref{CM} the following
\begin{Remark}
{\em Let $R$ be a bigraded Cohen-Macaulay $K$-algebra. Then
\[
H^{k}_P(\omega_R)\iso H^{d-k}_Q(R)^\vee \quad\text{for all $k$.}
\]}
\end{Remark}

Recall that for a finitely generated graded $R$-module $N$ one has
that $\dim_{R_0} N_j$ as well as $\depth_{R_0} N_j$ is constant
for large $j$, see \cite[Proposition, 2.5]{B}. In fact, if $N$ is
Cohen-Macaulay, then $\lim_{j\to \infty} \depth_{R_0} N_j=\dim
N-\dim N/P_{0}N$ as shown in \cite{HH}.  We call these constants
the limit depth and limit dimension, respectively. Using this fact
we have

\begin{Corollary}
\label{CMtame}
 Let $M$ be a bigraded Cohen-Macaulay $R$-module of
dimension $s$. We set $N=H^s_{R_+}(M)^\vee$, and put  $t_0=
\lim_{j\to \infty} \depth_{R_0} N_j$ and $s_0=\lim_{j\to
\infty}\dim_{R_0}N_j$. Then the $R$-modules $H^{j}_Q(M)$ are tame
for all $j\leq s-s_0$ and $j\geq s-t_0$.
\end{Corollary}
\begin{proof}
 We see that $H^{s-i}_P(N)_j=H^{s-i}_{P_0}(N_j)\neq 0$  for $j\gg 0$ if $i=s-s_0$  and $i=s-t_0$,  and
also $H^{s-i}_P(N)_j=H^{s-i}_{P_0}(N_j)= 0$ for $j\gg 0$ if
$i<s-s_0$ and $i>s-t_0$. Therefore  by Corollary \ref {CM} we have
the desired conclusion.
\end{proof}

\begin{Remark}
{\em In view of Corollary \ref{CM}  the   tameness conjecture of
Brodmann and Hellus would hold for bigraded  Cohen-Macaulay
modules provided one could show that for any bigraded
Cohen-Macaulay $R$-module $N$ and all integers $i$ there exists an
integer $j_0$ such that either
\[
H^i_{P_0}(N_j)=0  \quad \text{for all}\quad j\geq j_0, \quad
\text{or}\quad  H^i_{P_0}(N_j)  \neq 0\quad \text{for all}\quad
j\geq j_0.
\]
More generally one is tempted to conjecture the following: Let $M$
be a finitely generated bigraded $R$-module, and $W$ a finitely
generated $R_0$-module. Then for all $i$ there exists an integer
$j_0$ such that either
\[
\Ext^i_{P_0}(N_j,W)=0  \quad \text{for all}\quad j\geq j_0, \quad
\text{or}\quad  \Ext^i_{P_0}(N_j,W)  \neq 0\quad \text{for
all}\quad j\geq j_0.
\]
 A similar,  even more general conjecture can be made when $R$ is a
finitely generated positively graded $R_0$-algebra where $R_0$ is
Noetherian, $W$ is a finitely generated $R_0$-module and $N$ a
finitely generated graded $R$-module. In still another variation
of the conjecture one could replace $\Ext^i$ by $\Tor^i$.}
\end{Remark}

\begin{Corollary}
\label{structure1} Assume $R_0$ is Cohen-Macaulay and  $M$ is a
bigraded Cohen-Macaulay $R$-module of dimension $s$. We set
$N=H^s_{R_+}(M)^\vee$. Then
\begin{enumerate}
\item[(a)]  for all $k$ and $j$ we have the following isomorphism
of graded $R_0$-modules
\[\Ext_{R_0}^{d-k}(N_j,\omega_{R_0})\iso
H^{s-k}_Q(M)_{-j},
\]

where $d=\dim R_0$.

 \item[(b)] $\dim H^{s-k}_Q(M)_{-j}\leq k$ for all $k$ and
$j$.
\end{enumerate}
\end{Corollary}

\begin{proof} Corollary \ref{CM} implies that
$(H^{s-k}_Q(M)_{-j})^\vee\iso
(H^{s-k}_Q(M)^\vee)_j=H_{P_0}^k(N_j)$. Thus the local duality
theorem yields $H^{s-k}_Q(M)_{-j}\iso H_{P_0}^k(N_j)^\vee\iso
\Ext_{R_0}^{d-k}(N_j,\omega_{R_0})$, as desired.

Finally by \cite[Corollary,\ 3.5.11(c)]{BH} one has $\dim_{R_0}
\Ext_{R_0}^{d-k}(N_j,\omega_{R_0})\leq k$. This proves statement
(b).
\end{proof}

Let $N\neq 0$ be a graded $R_0$-module. We set  $a(N)=\inf\{i\:
N_i\neq 0\}$ and $b(N)=\sup\{i\: N_i\neq 0\}$. If $N=0$ we set
$a(N)=\infty$ and $b(N)=-\infty$.

Recall that the {\em regularity} of $N$ is defined to be
\[
\reg N=\max\{b(H_{P_0}^k(N))+k\: k=0,1,\ldots\}.
\]

With the assumptions and notation introduced in Corollary \ref{CM}
we therefore have
\[
\reg(N_j)=-\min\{a(H^{s-k}_Q(M)_{-j})-k\: k=0,1,\ldots\}.
\]
In \cite{CHT} and \cite{K} it is shown that $\reg(N_j)$ is bounded
above by  a linear function of $j$. Thus in view of the preceding
formula we get

\begin{Corollary}
\label{reg} Let $M$ be a Cohen-Macaulay $R$-module. Then there
exist integers $c$ and $d$ such that $a(H^{k}_Q(M)_{j})\geq cj+d$
for all $k$ and all $j$.
\end{Corollary}

If the dimension and the depth of $M$ differ at most by 1 or $\dim
R_0\leq 1$ one obtains

\begin{Proposition}
\label{tame}
The following statements hold:
\begin{enumerate}
\item[(a)] if $\dim M=s$ and $\depth M=s-1$, then  we obtain the
long exact sequence
\begin{eqnarray*}
\dots \rightarrow H^{m-j-2}_P(H^{s-1}_{R_+}(M)^\vee) \rightarrow H^{m-j}_P(H^s_{R_+}(M)^\vee)
\rightarrow  H^{s-m+j}_Q(M)^\vee \rightarrow \\ H^{m-j-1}_P(H^{s-1}_{R_+}(M)^\vee)\rightarrow
H^{m-j+1}_P(H^s_{R_+}(M)^\vee)\rightarrow \dots ;
\end{eqnarray*}

\item[(b)] if $\dim R_0=0$, then $H^{i}_{R_+}(M)\iso H^{i}_Q(M)$
for all  $i$;

\item[(c)] if $\dim R_0=1$, then for all $i$ we have the short exact sequence
 \[
0 \rightarrow H^{1}_P(H^{i+1}_{R_+}(M)^\vee) \rightarrow H^{i}_Q(M)^\vee \rightarrow H^{0}_P(H^{i}_{R_+}(M)^\vee)\rightarrow 0.
\]
\end{enumerate}

\end{Proposition}
\begin{proof}
We first prove (a). Our hypotheses imply the  following exact
sequences
\[
0 \rightarrow E^\infty_{s,j} \rightarrow H^{s-m+j}_Q(M)^\vee \rightarrow E^\infty_{s-1,j+1} \rightarrow 0,
\]
\[
0 \rightarrow E^\infty_{s-1,j+1}\rightarrow E^2_{s-1,j+1} \stackrel
{}\rightarrow  E^2_{s,j-1}\rightarrow
E^\infty_{s,j-1}\rightarrow 0.
\]
Putting these two exact sequences together we get the long exact sequence
\begin{eqnarray*}
\dots \rightarrow E^2_{s-1,j+2} \rightarrow E^2_{s,j}\rightarrow H^{s-m+j}_Q(M)^\vee \rightarrow E^2_{s-1,j+1}\rightarrow \\ E^2_{s,j-1}\rightarrow  H^{s-m-1+j}_Q(M)^\vee \rightarrow E^2_{s-1,j}\rightarrow  E^2_{s,j-2}\rightarrow \dots,
\end{eqnarray*}
which yields the desired exact sequence.

For the proof (b) we set $N=H^i_{R_+}(M)^\vee$. Since $\dim
R_0=0$, it follows that $N_k$ is a finitely generated $R_0$-module
of finite length. Thus $H^{m-j}_{P}(N)_k=H^{m-j}_{P_0}(N_k)=0$ for
all $k$ and $j<m$ and hence  $E^2_{i,j}=0$ for all $i$ and $j\neq
m$. Therefore we have $N=H^0_{P}(N)\iso H^i_Q(M)^\vee$ for all
$i$, and so $H^{i}_{R_+}(M)\iso H^{i}_Q(M)$ for all $i$. In order
to prove (c) we again set $N=H^i_{R_+}(M)^\vee$. Since $\dim
R_0=1$, it follows that $H^{m-j}_{P}(N)_k =H^{m-j}_{P_0}(N_k)=0$
for all $k$ and $j<m-1$ and hence $E^2_{i,j}=0$ for all $i$ and
$j\neq m, m-1$. Thus for all $i$ we get the exact sequence
\[
0 \rightarrow E^\infty_{i+1,m-1} \rightarrow H^{i}_Q(M)^\vee \rightarrow E^\infty_{i,m} \rightarrow 0.
\]
Since $E^\infty_{i+1,m-1}=E^2_{i+1,m-1}$ and $E^\infty_{i,m}=E^2_{i,m}$ for all $i$, the result follows.
\end{proof}

As a simple consequence of Proposition \ref{tame} (b),(c) we
obtain the following tameness result due to \cite[Theorem, 4.5]{Br}.

\begin{Corollary}
\label{tame2} Let $\dim R_0\leq 1$. Then $H^i_Q(M)$ is tame for
all $i$.
\end{Corollary}
\begin{proof}
First we assume that $\dim R_0=0$. Since any Artinian graded
$R$-module is tame, the result follows from Proposition \ref{tame}
(b).

Now we assume that $\dim R_0=1$. Let $N$ be a finitely generated
bigraded $R$-module. By Proposition \ref{tame} (c) it is enough to
prove that there exists an integer $j_0$ such that for $i=0,1$ one
has:
\[
H^i_{P_0}(N_j)=0   \quad \text{for all }\quad j\geq j_0,  \quad
\text{or}\quad  H^i_{P_0}(N_j)\neq 0 \quad \text{for all}\quad
j\geq j_0.
\]
We set
$t_0=
\lim_{j\to \infty} \depth_{R_0} N_j$ and $s_0=\lim_{j\to
\infty}\dim_{R_0}N_j$.
Then $H^0_{P_0}(N_j)\neq 0$  for $j\gg 0$ if $t_0=0$, and
$H^0_{P_0}(N_j)= 0$ for $j\gg 0$ if $t_0\neq 0$. Similarly,
$H^1_{P_0}(N_j)\neq 0$  for $j\gg 0$ if $s_0=1$, and
$H^1_{P_o}(N_j)$ for $j\gg 0$ if $s_0=0$.
\end{proof}

Finally we want to mention two standard 5-term exact sequences
arising from our spectral sequence

\begin{Proposition} There is a  5-term exact sequence for the
corner $(t,0)$

\begin{eqnarray*}
H^{t+2-m}_Q(M)^\vee \rightarrow H^{m-2}_P(H_{R_+}^t(M)^\vee) \rightarrow H^m_P(H_{R_+}^{d+1}(M)^\vee) \rightarrow\\
H^{t+1-m}_Q(M)^\vee \rightarrow H^{m-1}_P(H_{R_+}^t(M)^\vee) \rightarrow 0,
\end{eqnarray*}
and a 5-term exact sequence for the corner $(s,m)$
\begin{eqnarray*}
H^{s}_Q(M)^\vee \rightarrow H^{0}_P(H_{R_+}^{s-1}(M)^\vee)\rightarrow H^2_P(H_{R_+}^s(M)^\vee)\rightarrow\\
H^{s-1}_Q(M)^\vee \rightarrow H^{0}_P(H_{R_+}^{s-1}(M)^\vee)\rightarrow 0.
\end{eqnarray*}
\end{Proposition}

\end{document}